\documentclass[10pt]{elsarticle}

\def\RR{{\mathbb R}}
\def\NN{{\mathbb N}}

\def\CC{{\mathbb C}}

\def\beqns{\begin{eqnarray*}}
\def\eeqns{\end{eqnarray*}}
\def\beqn{\begin{eqnarray}}
\def\eeqn{\end{eqnarray}}

\def\proof{\noindent {\bf Proof:} \, }
\def\endproof{{\hfill $\blacksquare$}}
\def\mendproof{{\qquad \blacksquare}}

\def\erdos-1{Erd\H{o}s and Tur\'{a}n}

\newcommand{\sph}{\mathbb S}

\newcommand{\cH}{{\mathcal H}}
\newcommand{\homp}{\mathrm{hom}}
\newcommand{\harm}{\mathrm{harm}}

\newcommand{\cU}{{\mathcal U}}

\newcommand{\cP}{{\mathcal P}}

\newcommand{\ci}{{\rm i}}

\newtheorem{theorem}{Theorem}[section]
\newtheorem{lemma}[theorem]{Lemma}

\newtheorem{proposition}[theorem]{Proposition}
\newtheorem{corollary}[theorem]{Corollary}

\def\be{{\bf e}}
\def\bz{{\bf z}}
\def\bw{{\bf w}}
\def\bx{{\bf x}}
\def\by{{\bf y}}
\def\span{{\rm span} \,}
\def\dlangle{{\langle \! \langle}}
\def\drangle{{\rangle \! \rangle}}

\usepackage{amssymb}
\begin{document}

\begin{frontmatter}



\title{A Multiplier Version of the Bernstein Inequality on the Complex Sphere}


\author{Jeremy Levesley and Alex Kushpel }

\address{Department of Mathematics, University of Leicester,
University Road, Leicester, LE1 7RH}

\begin{abstract}

We prove a multiplier version of the Bernstein inequality on the complex sphere.
Included in this is a new result relating a bivariate sum involving Jacobi polynomials
and Gegenbauer polynomials, which relates the sum of reproducing kernels on spaces of
polynomials irreducibly invariant under the unitary group, with the reproducing kernel
of the sum of these spaces, which is irreducibly invariant under the
action of the orthogonal group.

\end{abstract}

\begin{keyword}

Bernstein Inequality, Complex Sphere, Multiplier, Jacobi Polynomial


\end{keyword}

\end{frontmatter}

% \linenumbers

\section{Introduction and preliminaries} \label{intro}

In this article we prove a multiplier version of the Bernstein
inequality of the type proved by Ditzian \cite{ditzian}. 
Since the restriction to a
geodesic of a polynomial on a complex sphere is just a trigonometric
polynomial on a circle, we immediately have a {\sl tangential
Bernstein inequality}
\[
\|D_u p_n \|_\infty \le n \|p_n\|_\infty,
\]
where $D_u$ is the tangential derivative in the direction of $u$ and $p_n$
is any polynomial of degree $n$. For more information on tangential Bernstein
inequalities on algebraic manifolds see e.g. Bos et. al. \cite{bos}.
An important stepping stone for this proof is Theorem~\ref{repronew}, in
which prove a new bivariate summation formula for Jacobi polynomials.

%In follow up papers we will use the Bernstein inequality proved here to
%prove results concerning simultaneous approximation on the complex sphere
%using positive definite zonal kernels, a Kolmogorov inequality, and $m$-term approximation.

We follow Koornwinder \cite{koornwinder} in our description of the
complex sphere, and the harmonic analysis thereof. Let $\CC^q$ be
$q$-dimensional complex space. We will denote vectors in $\CC^q$ by
$\bz = (z_1,z_2,\cdots,z_q)$. Let the inner product of two vectors $\bw, \bz \in \CC^q$ be
\[
\langle \bw,\bz \rangle = \sum_{j=1}^q w_j \overline{z}_j,
\]
and the length of a vector be $|\bz| = \langle \bz,\bz \rangle^{1/2}$. Let
\[
\sph^{2q} = \{ \bz \in \CC^q : | \bz | = 1 \},
\]
be the sphere in $\CC^q$. We note here that $\sph^{2q}$ has topological
dimension $2q-1$, but that we keep with the established notation so as not
to confuse the reader. Let $d(\bw,\bz)$ be the geodesic distance between $\bw$ and $\bz$ on $\sph^{2q}$.

The complex sphere is invariant under the action of the unitary
group $\cU_q$, the group of $q \times q$ complex matrices $U$ which
satisfy
\[
U U^* = I_q,
\]
where $U^*_{ij} = \overline{U_{ji}}$, $i,j=1,\cdots,q$.

Using the polar form for a complex number we can write $\bz \in \sph^{2q}$ in the form
\[
\bz = (r_1e^{\ci \phi_1},r_2e^{\ci \phi_2},\cdots,r_qe^{\ci \phi_q} ),
\]
where $\sum_{j=1}^q r_j^2 = 1$. If we set $r_1=\cos \theta$, we can write
\begin{equation}
\bz = \cos \theta e^{\ci \phi} \be_1 + \sin \theta \bz', \label{stform}
\end{equation}
where $\be_k$ is the unit vector in the $k$th coordinate, and $\bz' \in \sph^{2(q-1)}$.
Here $\phi=\phi_1$, (obviously) $\sin \theta = \sqrt{r_2^2+\cdots+r_q^2}$, and
\[
\bz'= (\sin \theta)^{-1} (r_2e^{\ci \phi_2},\cdots,r_qe^{\ci \phi_q} ).
\]

We can easily verify that $\sph^{2q} = \{ U \be_1, \; U \in \cU_q \}$.
Thus, for any $\bz \in \sph^{2q}$, there exists a $U \in \cU_q$ such that $U \be_1 = \bz$.
We call this action of $\cU_q$ on $\sph^{2q}$ {\em transitive}.  Now it is clear that
if we view $\cU_{q-2}$ as acting on the orthogonal complement of $\be_1$, then $\be_1$
remains fixed under this action. Thus we can write
\[
\sph^{2q} = {\cU_q \over \cU_{q-1}}.
\]
On the real sphere we are accustomed to the idea that the polynomials on the
sphere may be orthogonally decomposed into subspaces of spherical harmonics,
each of which is invariant under the action of the orthogonal group.
For the complex sphere the picture is not so straightforward.  Now we wish to
identify the spaces of polynomials which are minimally invariant under the
action of the unitary group, and this issue is discussed in Section~\ref{harmonic}.

Let $d \mu_{2q}$ be the $\cU_q$-invariant normalised measure on the sphere, and define the
inner product of $f,g$, two functions on $S^{2q}$, by
\[
\dlangle f,g \drangle = \int_{\sph^{2q}} f \overline g d\mu_{2q}.
\]
Let us define the family of $L_r$ norms on $\sph^{2q}$:
\[
\| f \|_r = \left \{ \begin{array}{ll}
\left ( \int_{\sph^{2q}}| f |^r d \mu_{2q} \right )^{1/r}, & 1<r<\infty, \\
{\rm ess \; sup} | f |, & r=\infty.
\end{array} \right.
\]
In this paper we will be discussing $\cU_q$ invariant kernels on
$\sph^{2q}$. These are kernels $\kappa:\sph^{2q} \times \sph^{2q}
\rightarrow \CC$, such that $\kappa(U \bx,U\by)=\kappa(\bx,\by)$ for
all $U \in \cU_q$. Previous results of Ditzian \cite{ditzian} have
been valid for two-point homogeneous spaces. These are spaces which
for pairs of points which are equidistant, there is a single
isometry which maps one pair to the other (see Wang \cite{wang} for more information).  For two points spaces,
the geodesic distance is a function of the inner product in the
ambient space. A consequence of this is that all isometrically
invariant kernels are univariate functions of distance.

For the complex spheres this is not the case. However, we does have
the following analogous property. Suppose we have pairs of points
$\bx_1,\by_1$ and $\bx_2,\by_2$, with $\langle \bx_1,\by_1 \rangle =
\langle \bx_2,\by_2 \rangle$.  Since the unitary group acts
transitively on the complex sphere, there exist $U_1,U_2 \in \cU_q$
such that $U_1 \bx_1 = U_2 \bx_2 = \be_1$. Recalling (\ref{stform}),
and using the fact that $\cU_{q-1}$ acts transitively on
$\sph^{2q-2}$, we know there exists $U^\prime \in \cU_{q}$, such
that $U^\prime U_1 \by_1 = U_2 \by_2$, and $U^\prime \be_1 = \be_1$.
Hence, $U_2^{-1} U^\prime U_1 \bx_1 =  \bx_2$, and $U_2^{-1}
U^\prime U_1 \by_1 =  \by_2$. Hence, we conclude that if $\langle
\bx_1,\by_1 \rangle = \langle \bx_2,\by_2 \rangle$ there exists $U
\in \cU_{q}$ such that $U \bx_1=\bx_2$ and $U \by_1=\by_2$. This is
analogous to the two point homogeneous property of the reals
spheres. A straightforward consequence of this is that if $\kappa$
is $\cU_q$ invariant $\kappa(\bx_1,\bx_2)=\kappa(U \bx_1,U
\bx_2)=\kappa(\by_1,\by_2)$, so that $\kappa$ is invariant on points
with $\langle \bx_1,\by_1 \rangle = \langle \bx_2,\by_2 \rangle$.
Thus we have

\begin{lemma} \label{invar}
If $\kappa$ is a $\cU_q$-invariant kernel then
\[
\kappa(\bx,\by)=\psi(\langle \bx,\by \rangle)
\]
for some univariate function $\psi$.
\end{lemma}

We can define a convolution of an arbitrary $f \in L_1(\sph^{2q})$
function, with a $\cU_q$-invariant kernel $h \in L_1(\sph^{2q})$ function:
\[
f \ast \kappa (\bx) = \int_{\sph^{2q}} f(\by) \psi(\langle \bx,\by \rangle)d\mu(\by).
\]
It is observed in \cite{koornwinder} that we may view $\CC^q$ with
typical point
\[
\bz=(x_1+\ci y_1, x_2+\ci y_2,\cdots,x_q+\ci y_q)
\]
as a $2q$-dimensional real space with variables
\[
\bw=(x_1,y_1,x_2,y_2,\cdots,x_q,y_q).
\]
The inner product of two vectors $\bw$ and $\bw^{\prime}$ in this space is
\[
(\bw,\bw^\prime)=\sum_{j=1}^q (x_j x_j^\prime+y_j y_j^\prime) = \Re \langle \bz,\bz^{\prime} \rangle,
\]
with $\bz^{\prime} = (x_1^\prime+\ci y_1^\prime, x_2^\prime+\ci y_2^\prime,\cdots,x_q^\prime+\ci y_q^\prime)$.
Hence, a point with standard representation (\ref{stform}) on the complex sphere, has geodesic distance
$\cos^{-1} (\cos \theta \cos \phi)$ from the north pole on the associated real sphere.

%%%%%%%%%%%%%%%%%%%%%%%%%%%%%%%%%%%%%%%%%%%%%%%%%%%%%%%%%%%%%%%%%%%%%%%%%%

%%%%%%%%%%%%%%%%%%%%%%%%%%%%% 2 %%%%%%%%%%%%%%%%%%%%%%%%%%%%%%%%%%%%%%%%%%

\section{Harmonic analysis} \label{harmonic}
To start with it might be informative to briefly discuss harmonic analysis on the circle
as a subset of the complex numbers as opposed to a subset of $\RR^2$. In the former case
we complex Fourier series with a basis $\{1, z^k,\overline{z}^k\}$, $k=1,\cdots$.
The unitary group in this case is just the unit circle in the complex numbers.
Invariant subspaces under the action of the unitary group are just the one dimensional
spaces, constants, $\span \{ z^k \}, \span \{ \overline{z}^k \}$, $k=1,2,\cdots$.
For the latter case we have a basis $\{1, \Re(z^k), \Im(z^k)\}$, $k=1,2,\cdots$.
The subspaces which are invariant under $2 \times 2$ orthogonal matrices are,
constants, $\span \{ \Re(z^k), \Im(z^k) \}$, $k=1,2,\cdots$, which are two dimensional.
Hence we see that the use of the unitary matrices, as opposed to the orthogonal
matrices has given us a finer division of the polynomial spaces.

In this spirit let us define the space $\cP(m,n)$ of homogeneous polynomials in
$\CC^q$ as those of the form
$P(\bz,\overline{\bz})=P(z_1,z_2,\cdots,z_q,\overline{z}_1,\overline{z}_2,\cdots,\overline{z}_q)$, satisfying
\[
P(\alpha \bz, \overline{\beta \bz}) = \alpha^m \overline{\beta}^n P(\bz,\overline{\bz}), \quad m,n \in \NN, \quad \alpha,\beta \in \CC.
\]
Here we are regarding $\bz$ and $\overline{\bz}$ formally as different variables,
though this is not really the case.
Then we define $\homp(m,n)$, the space of homogeneous polynomials on the sphere,
to be the restriction of $\cP(n,m)$ to the sphere via
\[
p(\bz) = P(\bz,\overline{\bz}), \quad \bz \in \sph^{2q}.
\]
Since on the sphere $z_1 \overline{z}_1+z_1
\overline{z}_2+\cdots+z_q \overline{z}_q=1$, we have $\homp(m-1,n-1)
\subset \homp(m,n)$. We define the space of harmonic polynomials
$\harm(m,n)=\homp(m,n) \cap \homp(m-1,n-1)^\perp$, where
orthogonality is with respect to the inner product above. In
\cite{koornwinder} these polynomials are defined in terms of the
Laplace operator, but Theorem 3.4 therein:
\[
\homp(m,n) = \oplus_{k=0}^{\min(m,n)} \harm(m-k,n-k),
\]
tells us that our definition is equivalent. For ease of notation
let us write $m \wedge n = \min(m,n)$.

From \cite{koornwinder} we know that the dimension of $\harm(m,n)$ is
\begin{equation}
d_{m,n}={(m+n+q-1) (m+q-2)! (n+q-2)! \over m! n! (q-1)! (q-2)! }.
\end{equation}
Now, let $\cH_l = \oplus_{k=0}^{l} \harm(l-k,k)$ be the harmonic space of degree $l$.
We can compute the dimension $d_l$ of $\cH_l$ directly by summation, but also we have
that it has the same dimension as the space of spherical harmonics in $\RR^{2q}$,
which from e.g. M\"uller \cite[Page 4]{muller} is
\begin{equation} \label{dimhom}
d_l = {l+2q-1 \choose l} - {l+2q-3 \choose l} =  {2 \over (2q-2)!} { (l+q-1) (l+2q-3)! \over l!}.
\end{equation}
The dimension of the full polynomial space $\cP_n = \oplus_{l=0}^n \cH_l$ is
\begin{equation} \label{tndef}
t_n = {1 \over (2q-1)!} { (2n+2q-1) (n+2q-2)! \over
 n!}.
\end{equation}

Let $k_1,k_2,\cdots,k_{d_{m,n}}$ be an orthonormal basis for $\homp(m,n)$.
Then the reproducing kernel for projection onto $\homp(m,n)$ is
 \begin{equation} \label{rep}
 \kappa_{m,n}(\bz,\bw) = \sum_{j=1}^{d_{m,n}} k_j(\bz) \overline k_j(\bw).
 \end{equation}
 It is straightforward to show that this kernel is $\cU_q$-invariant.
 Similarly the reproducing kernel for $\cH_l$,
 \[
 h_{l}(\bz,\bw) = \sum_{k=0}^{l} \kappa_{l-k,k}(\bz,\bw),
 \]
 is $\cU_q$-invariant.

Since the reproducing kernels are $\cU_q$-invariant we have, from Lemma~\ref{invar}, that
\[
\kappa_{m,n}(\bz,\bw)=\psi(\langle \bz,\bw \rangle ),
\]
for some univariate function $\psi$. In order to determine $\psi$ we need
to use its orthogonality properties.

In terms of the standard representation (\ref{stform}) we can write the
surface element on $\sph^{2q}$ as
\[
d \mu_{2q} ({\bf z}) = \cos \theta (\sin \theta)^{2q-3} d \theta d\phi d \mu_{2q-2},
\]
since surface area on $r \mathbb{S}^{2q-2}$ scales like $r^{2q-3}$.
If we make the change of variable $t=\cos \theta$ then we see that the measure
$\cos \theta (\sin \theta)^{2q-3} d \theta = t(1-t^2)^{q-2} dt$ arises.
Thus we might expect the reproducing kernels for the harmonic subspaces,
which are orthogonal, to be related to orthogonal polynomials with a
weight $(1-t^2)^{q-2}$, and indeed this is the case.

From \cite{koornwinder} we have the following representation of the
reproducing  kernels for the irreducible polynomial spaces
$\cH(m,n)$,
\begin{equation} \label{repromn}
k_{m,n}(\bx,\by)=d_{m,n} e^{\ci(m-n)\phi} (\cos
\theta)^{|m-n|} { P_{m \wedge n}^{(q-2,|m-n|)}(\cos (2 \theta) )
\over P_{m \wedge n}^{(q-2,|m-n|)}(1 ) },
\end{equation}
where $\langle \bx,\by \rangle = \cos \theta e^{\ci \phi}$.
Here $P_{j}^{(\alpha,\beta)}$ is the degree $j$ Jacobi polynomial
which is orthogonal with respect to the weight $(1-t)^\alpha (1+t)^\beta$.
For ease of notation we will now write $\kappa_l(\theta,\phi)$ instead of $\kappa_l(\bx,\by)$

As stated in the introduction, we can also view the complex sphere as a real sphere.
The harmonics in $\homp(m,n)$ are complex harmonics on the real sphere of degree
$m+n$; see \cite{koornwinder}. The associated real sphere is of dimension $2q-1$.
Hence, we have the following reproducing kernel formula
\begin{equation} \label{rep2}
\sum_{m+n=l} \kappa_{m,n}(\bz,\bw) = h_l(\bz,\bw)=d_l { P_l^{(q-1)} (( \bz,\bw )) \over P_l^{(q-1)} (1)},
\end{equation}
where $P_l^{(\sigma)}, l \ge 0$ are the Gegenbauer polynomials which
are orthogonal with respect to the weight $(1-t^2)^{\sigma-1/2}$.
We normalise the Gegenbauer polynomials by
\[
P_l^{(\sigma)}(1) = {l+2 \sigma-1 \choose l}.
\]
Here we interpret $\bz$ and $\bw$ as points on the real sphere, and if
$\langle \bz,\bw \rangle = e^{\ci \phi} \cos \theta$ then $( \bz , \bw) = \cos \theta \cos \phi$
(see the closing remarks of Section~\ref{intro}).

In mind of (\ref{rep}) and (\ref{repromn}), we have the following interesting
(and we believe new) formula relating Jacobi and Gegenbauer polynomials.

\begin{theorem} \label{repronew}
For $d \ge 1$ and $l \ge 0$
\[
\sum_{m+n=l} d_{m,n} e^{\ci(m-n)\phi} (\cos
\theta)^{|m-n|} { P_{m \wedge n}^{(q-2,|m-n|)}(\cos (2 \theta) )
\over P_{m \wedge n}^{(q-2,|m-n|)}(1 ) } = d_l { P_l^{(q-1)} (\cos \theta \cos \phi) \over P_l^{(q-1)} (1)}.
\]
\end{theorem}

We wish to define multiplier (pseudodifferential) operators via their action on the
harmonic subspaces $\cH_l$. Let $M_l$ be the orthogonal projector from $L_2(\sph^{2q}) \rightarrow \cH_l$,
$l=0,1,\cdots$. The kernel of this projection is $h_l$, so that
\[
M_l f = f \ast h_l.
\]
Let $\lambda_l, l=0,1,\cdots,$ be a sequence of increasing real numbers.
Then, for $f \in L_1(\sph^{2q})$ (which thus has a formal Fourier expansion),
the multiplier operator
\[
\Lambda f = \sum_{l=0}^\infty \lambda_l M_l f.
\]
In Theorem~\ref{bernstein}, in Section~\ref{bern}, we will show that for $p \in \cP_n$,
\[
\| \Lambda p \|_r \le \lambda_n \| p \|_r, \quad 1 \le r \le \infty.
\]

%%%%%%%%%%%%%%%%%%%%%%%%%%%%%%%%%%%%%%%%%%%%%%%%%%%%%%%%%%%%%%%%%%%%%%%%%%

\section{Cesaro means for reproducing kernels} \label{cesaro}

In order to prove Theorem~\ref{bernstein} we observe that for $p \in \cP_n$,
\begin{eqnarray*}
\Lambda p & = & K_m \ast p,  \quad m \ge n,
\end{eqnarray*}
where
\begin{equation} \label{km}
K_m = \sum_{l=0}^n \lambda_l h_l + \sum_{l=n+1}^m \tilde \lambda_l h_l,
\end{equation}
where the numbers $\tilde \lambda_l$, $l=n+1,\cdots,m$ are available for us to choose.
Let us define the sequence $\rho_l=\lambda_l$, $l=0,1,\cdots,n$ and
$\rho_l=\tilde \lambda_l$, $l=n+1,\cdots,m$.

Using Young's inequality
\[
\| f \ast g \|_p = \| f \|_p \| \| g \|_1,
\]
we are directed towards the computation of the 1-norms of the kernels $K_m$,
which we achieve via the Cesaro means of $h_l$
\[
S_m^\delta   = {1 \over C_m^\delta} \sum_{l=0}^m C_{m-l}^\delta h_l,
\]
where
\[
C_k^\delta = {k+\delta \choose k}  \asymp k^\delta, \quad k=0,1,\cdots,m.
\]
Before we proceed we need a preliminary technical lemma:
\begin{lemma} \label{dmmk}
Let (see \cite[4.7.15]{szego})
\begin{eqnarray*}
\lefteqn{\gamma_l^{(\sigma)} = \int_{-1}^1 (1-t)^{\sigma-1/2} |P_l^{(\sigma)}(t) |^2 dt} \\
& = & {2^{1-2\sigma} \pi \over  (\Gamma(\sigma))^2} {\Gamma(l+2\sigma) \over (l+\sigma) \, l! }.
\end{eqnarray*}
Then, for $l > 0$,
\[
d_{l} = \dim \cH_l = { 2^{(3-2q)} \pi (2q-3)!\over (q-1)! (q-2)!} {(P_l^{(q-1)}(1))^2 \over \gamma_l^{(q-1)}}.
\]
\end{lemma}
\proof Starting from (\ref{dimhom}), a straightforward calculation gives us
\begin{eqnarray*}
\lefteqn{d_{l}={2 \over (2q-2)!} { (l+q-1) (l+2q-3)! \over l!} }\\
 & = & {2 \over (2q-2)!} { 2^{(3-2q)} \pi ((2q-3)!)^2 \over ((q-2)!)^2} \left
 \{  { ((q-2)!)^2 \over 2^{(3-2q)} \pi} { (l+q-1) l! \over  (l+2q-3)! } \right \}
 \left ( {(l+2q-3)! \over l! (2q-3)!} \right )^2 \\
 & = & { 2^{(3-2q)} \pi (2q-3)!\over (q-1)!(q-2)!} {(P_l^{(q-1)}(1))^2 \over h_l^{(q-1)}}. \mendproof
\end{eqnarray*}

Using this last result and (\ref{rep2}) we see that
\[
S_n^\delta (\cos \psi)  = { 2^{(3-2q)} \pi (2q-3)!\over (q-1)! (q-2)!} {1 \over C_n^\delta}
\sum_{l=0}^n C_{l-n}^\delta  {P_l^{(q-1)}(1) P_l^{(q-1)}(\cos \psi) \over h_l^{(q-1)}},
\]
where $\cos \psi = \cos \theta \cos \phi$, in other words are essentially the
Cesaro means of the Gegenbauer polynomials.

Using Equation \cite[4.5.3]{szego} and Lemma~\ref{dmmk} we have the
following corollary of Theorem~\ref{repronew}:
\begin{corollary}
The reproducing kernel for $\cP_n$ is
\begin{eqnarray*}
r_n((\bz,\bw)) =\sum_{l=0}^n h_l & = & t_n {P_n^{(q+1/2,q-1/2)} (\cos \theta \cos \phi) \over P_n^{(q+1/2,q-1/2)}(1)},
\end{eqnarray*}
with $(\bz,\bw)  = \cos \theta \cos \phi$, where we recall that $t_n = \dim (\cP_n)$.
\end{corollary}

To estimate these we use the the following results which are given in
Bonami and Clerk \cite[Page 230]{bonami}.

\begin{proposition} \label{gegbnd}
If $0 \le \delta \le q$ then there is a constant $C$ such that,
\begin{equation} \label{bon}
S_n^\delta (\cos \psi) \le C \left \{ \begin{array}{ll}
n^{q-\delta-1} \psi^{-(q+\delta)}, & 3/n \le \psi \le \pi/4, \\
n^{2 q-1}, & 0 \le \psi \le 3/n.
\end{array} \right.
\end{equation}
\end{proposition}
In the remainder of this paper the number $C$ will be used to denote a
constant which is independent of $n$.

The main result of this section is
\begin{theorem} \label{thces}
For $0 \le \delta \le q$,
\[
\| S_n^\delta  \|_1 \le C \left \{ \begin{array}{ll}
n^{q-1-\delta}, & \delta \le q-2, \\
(\log n)^2, & \delta=q-1, \\
1, & \delta \ge q.
\end{array} \right.
\]
\end{theorem}

\proof We will provide a bound for
\begin{eqnarray*}
\| S_n^\delta \|_1 & = &  \int_{\sph^{2q}} | S_n^\delta (( \bz, \be ))| d \mu(\bz) \\
 & = & \int_{0}^{2 \pi} \int_0^{\pi/2} \cos \theta (\sin \theta)^{2q-3} | S_n^\delta (\cos \theta \cos \phi)| d \theta d \phi \\
 & = & 2 \int_{0}^{\pi} \int_0^{\pi/2} \cos \theta (\sin \theta)^{2q-3} | S_n^\delta (\cos \theta \cos \phi)| d \theta d \phi.
\end{eqnarray*}
Suppose that $Q$ is the region $\pi/4 \le \theta \le \pi/2$ or $\pi/4 \le \phi \le \pi$.
Then $\cos \theta \cos \phi \le 1/\sqrt{2}$. Setting $\cos \psi = \cos \theta \cos \phi$,
we have $\pi/4 \le \psi \le 3 \pi/4$. Since, from Proposition~\ref{gegbnd}, $S_\nu^\delta (\cos \psi)$
is bounded above for $\pi/4 \le \psi \le 3 \pi/4$ we have
\begin{equation} \label{i1}
I_1 = \int_{Q} \cos \theta (\sin \theta)^{2q-3} | S_n^\delta (\cos \theta \cos \phi)| d \theta d \phi \le C.
\end{equation}
We break the remaining integral into 4 parts, $(\theta,\phi) \in [0,1/n]^2$, $[0,1/n] \times [1/n,\pi/4]$,
$[1/n,\pi/4] \times [0,1/n]$ and $[1/n,\pi/4] \times [1/n,\pi/4]$ which we call $I_2,I_3,I_4$ and
 $I_5$ respectively. Firstly, since on $[0,1/n]^2$, $\cos \theta \cos \phi \ge \cos^2(1/n) \ge \cos (3/n)$
(this is easy to check), we have
\begin{eqnarray} \label{i2}
I_2 & = & \int_{0}^{1/n} \int_0^{1/n} \cos \theta (\sin \theta)^{2q-3} | S_n^\delta (\cos \theta \cos \phi)| d \theta d \phi \nonumber \\
& \le & C n^{2q-1} \int_0^{1/n} \left ( \int_0^{1/n} \cos \theta (\sin \theta)^{2q-3} d \theta \right ) d \phi \nonumber \\
& \le & C.
\end{eqnarray}
Now, for the remaining integrals we observe that
\begin{eqnarray*}
\psi & = & \arccos (\cos \theta \cos \phi) \\
& \le & C(1-\cos \theta \cos \phi)^{1/2},
\end{eqnarray*}
since $\arccos(z) \le C(1-z)^{1/2}$ as $z \rightarrow 1$. Now, for $0 \le \theta,\phi \le 1$,
\begin{eqnarray*}
1-\cos \theta \cos \phi & \ge & 1-\left (1-{\theta^2 \over 2} \right )\left (1-{\phi^2 \over 2} \right ) \nonumber \\
& = & {\theta^2 + \phi^2 \over 2}-{\theta^2 \phi^2\over 4} \nonumber \\
& \ge & {\theta^2 + \phi^2 \over 4}.
\end{eqnarray*}
If we use this last equation in (\ref{bon}), we see that
\begin{equation} \label{snd}
|S_n^\delta (\cos \psi)| \le C n^{q-\delta-1} (1-\cos \theta \cos \phi)^{-(q+\delta)/2} \le C n^{q-\delta-1} (\theta^2 + \phi^2)^{-(q+\delta)/2}.
\end{equation}
We note that for any $\alpha,\beta > 0$,
$$
(\theta^2 + \phi^2)^{-(\alpha+\beta)} \le \theta^{-2 \alpha} \phi^{-2 \beta}.
$$
We have
\begin{eqnarray*} \label{i3}
I_3 & = & 2 \int_{1/n}^{\pi/4} \int_{0}^{1/n} \cos \theta (\sin \theta)^{2q-3} | S_n^\delta (\cos \theta \cos \phi)| d \theta d \phi \\
& \le & C n^{q-1-\delta} \int_{0}^{1/n} \int_{1/n}^{\pi/4} \theta^{2q-3} (\theta^2 + \phi^2)^{-(q+\delta)/2} d \theta d \phi.
\end{eqnarray*}
Now if $\delta \le q-2$, then
\begin{eqnarray*}
\int_{1/n}^{\pi/4} \int_{0}^{1/n} \theta^{2q-3} (\theta^2 + \phi^2)^{-(q+\delta)/2} d \theta d \phi
& \le & \int_{1/n}^{\pi/4} \int_{0}^{1/n} (\theta^2 + \phi^2)^{q-3/2-(q+\delta)/2} d \theta d \phi \\
& \le & \int_{1/n}^{\pi/2} r^{(q-\delta-3)} r dr < C,
\end{eqnarray*}
using a change to polar coordinates. If $\delta=q-1$, we have
(assuming $q >1$ and using $(\theta^2 + \phi^2)^{-q+1} \le
\theta^{-2q+5/2} \phi^{-1/2}$)
\begin{eqnarray*}
\int_{0}^{1/n} \int_{1/n}^{\pi/4} \theta^{2q-3} (\theta^2 + \phi^2)^{-q+1/2} d \theta d \phi
& \le & \int_{0}^{1/n} \phi^{-1/2} d \phi \int_{1/n}^{\pi/4} \theta^{-1/2}  d \theta \\
& \le & C.
\end{eqnarray*}
If $\delta=q$, we have (using $(\theta^2 + \phi^2)^{-q+1} \le
\theta^{-2q+5/2} \phi^{-1/2}$)
\begin{eqnarray*}
\int_{0}^{1/n} \int_{1/n}^{\pi/4} \theta^{2q-3} (\theta^2 + \phi^2)^{-q} d \theta d \phi  & \le & \int_{0}^{1/n} \phi^{-1/2} d \phi \int_{1/n}^{\pi/4} \theta^{-5/2}  d \theta \\
& \le & C n.
\end{eqnarray*}
Putting these estimates for the integral into (\ref{i2}) we have
 \begin{equation} \label{i32}
I_3 \le C \left \{ \begin{array}{ll}
n^{q-1-\delta}, & \delta \le q-2, \\
1, & \delta \ge q-2. \\
\end{array} \right.
\end{equation}
More straightforwardly,
\begin{eqnarray*}
I_4 & = & 2 \int_{1/n}^{\pi/4} \int_{0}^{1/n} \cos \theta (\sin \theta)^{2q-3} | S_n^\delta (\cos \theta \cos \phi)| d \theta d \phi \\
& \le & C n^{q-1-\delta} \int_{1/n}^{\pi/4} \int_{0}^{1/n} \theta^{2q-3} (\theta^2 + \phi^2)^{-(q+\delta)/2} d \theta d \phi \\
& \le & C n^{q-1-\delta} \int_{1/n}^{\pi/4} \phi^{-q-\delta} d \phi \int_{0}^{1/n} \theta^{2q-3} d \theta \\
& \le & C n^{q-1-\delta+q+\delta-1-2q+2} = C.
\end{eqnarray*}
For the last integral
\begin{eqnarray} \label{i5}
I_5 & = & 2 \int_{1/n}^{\pi/4} \int_{1/n}^{\pi/4} \cos \theta (\sin \theta)^{2q-3} | S_n^\delta (\cos \theta \cos \phi)| d \theta d \phi \nonumber \\
& \le & C n^{q-1-\delta} \int_{1/n}^{\pi/4} \int_{1/n}^{\pi/4} \theta^{2q-3} (\theta^2 + \phi^2)^{-(q+\delta)/2} d \theta d \phi.
\end{eqnarray}
Now if $\delta \le q-2$, then
\begin{eqnarray*}
\int_{1/n}^{\pi/4} \int_{1/n}^{\pi/4} \theta^{2q-3} (\theta^2 + \phi^2)^{-(q+\delta)/2} d \theta d \phi
& \le & \int_{1/n}^{\pi/4} \int_{1/n}^{\pi/4} (\theta^2 + \phi^2)^{q-3/2-(q+\delta)/2} d \theta d \phi \\
& \le & \int_{1/n}^{\pi/2} r^{(q-\delta-3)} r dr < C.
\end{eqnarray*}
If $\delta=q-1$ we have
\begin{eqnarray*}
\int_{1/n}^{\pi/4} \int_{1/n}^{\pi/4} \theta^{2q-3} (\theta^2 + \phi^2)^{-(q+\delta)/2} d \theta d \phi
& \le & \int_{1/n}^{\pi/4}  \int_{1/n}^{\pi/4} \theta^{2q-3} (\theta^2 + \phi^2)^{-(q-1/2)} d \theta d \phi \\
& \le & \int_{1/n}^{\pi/4} \phi^{-1} d \phi \int_{1/n}^{\pi/4} \theta^{-1}  d \theta \\
& \le & C (\log n)^2.
\end{eqnarray*}
For the last case, $\delta=q$, we have
\begin{eqnarray*}
\int_{1/n}^{\pi/4} \int_{1/n}^{\pi/4} \theta^{2q-3} (\theta^2 + \phi^2)^{-(q+\delta)/2} d \theta d \phi
& \le & \int_{1/n}^{\pi/4}  \int_{1/n}^{\pi/4} \theta^{2q-3} (\theta^2 + \phi^2)^{-2q} d \theta d \phi \\
& \le & \int_{1/n}^{\pi/4} \phi^{-3/2} d \phi \int_{1/n}^{\pi/4} \theta^{-3/2}  d \theta \\
& \le & C n.
\end{eqnarray*}
Hence, substituting the above estimates into (\ref{i5}), we see that
\begin{equation} \label{i52}
I_5 \le C \left \{ \begin{array}{ll}
n^{q-1-\delta}, & \delta \le q-2, \\
(\log n)^2, & \delta=q-1, \\
1, & \delta=q.
\end{array} \right.
\end{equation}
A simple inspection of the bounds (\ref{i1}) to (\ref{i52}) tells us
that the bound for (\ref{i5}) is the largest, giving the required
result. Estimates for $I_6$ and $I_7$ can be obtained similarly. \endproof

%%%%%%%%%%%%%%%%%%%%%%%%%%%%%%%%%%%%%%%%%%%%%%%%%%%%%%%%%%%%%%%%%%%%%%%%

\section{The Bernstein inequality} \label{bern}

The Laplace Beltrami operator on the complex sphere has eigenspaces $\cH_l$, with
eigenvalue $\lambda_l=l(l+2q-1)$, $l=0,1,\cdots$; see e.g. \cite{tao}. Thus the
fractional order differential operator $\Lambda$ has multiplers $\lambda_l=(l(l+2q-1))^{\gamma/2}$.

Performing Abel summation $q+1$ times on (\ref{km}) we get
\begin{equation} \label{abel}
K_m = \sum_{l=1}^m \Delta^{q+1} \rho_l C_l^q S_l^q+\sum_{l=0}^q \Delta^{l} \rho_{m-l} C_{m-l}^l S_{m-l}^l,
\end{equation}
where $\Delta^0 \rho_k=\rho_k$, $\Delta \rho_k = \rho_k-\rho_{k+1}$, and $\Delta^j \rho_k =  \Delta (\Delta^{j-1} \rho_k)$, $j=2,3,\cdots$.

Now let us define
\[
g(x)=\left \{ \begin{array}{ll}
1, & 0 \le x < n, \\
1-C_{n,q} \int_n^x \left|(y-n)(2n-y)\right|^{q+1} dy, & n \le x \le 2n, \\
0, & x>2n,
\end{array} \right.
\]
where
\[
C_{n,q} = \left (\int_n^{2n} \left | (y-n)(2n-y) \right |^{q+1} dy \right )^{-1}
= {(2q+3)! \over n^{2q+3} ((q+1)!)^2}.
\]
This last equations follows by making the change of variable $ns=y-n$ in the above integral, giving
\begin{eqnarray*}
\int_n^{2n} \left | (y-n)(2n-y) \right |^{q+1} dy & = & n^{2q+3} \int_0^{1} s^{q+1}(1-s)^{q+1} ds \\
& = & n^{2q+3} B(q+2,q+2),
\end{eqnarray*}
where $B$ is the Beta function. We now use $B(u,v) = \Gamma(u)
\Gamma(v)/\Gamma(u+v)$; (see \cite[Page 258]{stegun}). We have $g \in {\rm
C}^{(q+1)}(\RR_+)$, and a simple computation shows that, for $1 \le
j \le q+1$,
\[
g^{(j)}(x) =
C_{q,n} \sum_{k=1}^{\left\lfloor (j+1)/2 \right\rfloor} \nu_k ((x-n)(2n-x))^{q+1-j+k} (3-2x)^{j+1-2k}, \quad n \le x \le 2n,
\]
and is zero otherwise. Here, the $\nu_k$ depend on $q$, and 
$\left\lfloor  \cdot \right\rfloor$ denotes the integer part. Hence, for $1 \le j \le q+1$, we can bound
\begin{eqnarray*}
| g^{(j)}(x) | & \le & C n^{-j}.
\end{eqnarray*}
Let $h(x)=(x(x+2q-1))^{\gamma/2}$, and $f=gh$, and $\rho_k=f(k)$. We
have then, for $0 \le k \le n$, $\rho_k=\lambda_k$. We observe that
$\rho \in {\rm C}^{(q+1)}(\RR_+)$, where $\RR_+ = \{ x : x \ge 0\}$.
We can estimate the difference
\[
|\Delta^l \rho_k| \le C \max_{x \in [k,k+l]} |f^{(l)}(x)|, \quad l=0,1,\cdots,q+1.
\]
Using Leibnitz rule we have, for $n \le x \le 2n+l$,
\begin{eqnarray*}
|f^{(l)}(x)| & = & \left |\sum_{j=0}^l {l \choose j} g^{(j)} h^{(l-j)}(x) \right | \\
& \le & C \sum_{j=0}^l n^{-j} n^{\gamma-l+j} \\
& \le & C n^{\gamma-l}.
\end{eqnarray*}
Therefore, for $0 \le l \le q+1$,
\[
|\Delta^l \rho_k| \le C \left \{
\begin{array}{ll}
(k+l)^{\gamma-l}, & 0 \le k \le n-l, \\
n^{\gamma-l+j}, & n-l \le k \le n+l, \\
0, & {\rm otherwise}.
\end{array} \right.
\]
If we substitute these estimates, with Theorem~\ref{thces} into
(\ref{abel}) we see that, if $m \ge 2n+q+1$ (and so $|\Delta^{l}
\rho_{m-l}| = 0$, for $l=0,1,\cdots,q$)
\begin{eqnarray*}
\| K_m \|_1 & \le & \sum_{l=1}^m |\Delta^{q+1} \rho_l| C_l^q \| S_l^q \|_1  \\
& \le & C \left ( \sum_{l=1}^{n-q-1} l^{\gamma-q-1} l^q  + \sum_{l=n-q}^{2n+q+1} n^{\gamma-q-1} l^q  \right ) \\
& \le & C n^\gamma.
\end{eqnarray*}

Thus we have
\begin{theorem} \label{bernstein}
For $\lambda_l=(l(l+2d-1))^{\gamma/2}$,
$$
\| \Lambda p \|_r \le C n^\gamma \| p \|_r, \quad 1 \le r \le \infty,
$$
for every $p \in \cP_n$.
\end{theorem}

%%%%%%%%%%%%%%%%%%%%%%%%%%%%%%%%%%%%%%%%%%%%%%%%%%%%%%%%%%%%%%%%%%%%%%%%%%



%% References


%\bibliographystyle{model1c-num-names}

\begin{thebibliography}{30}


%\bibitem{25}

%%%%%%%%%%%%%%%%%%%%%%%%%%%%%%%%%%%%%%%%%%%%%%%%%%%%%%%%%%%%%%%%%%%%%%

\bibitem{stegun}  M. Abramowitz and I. A. Stegun, eds. Handbook of Mathematical Functions with Formulas, Graphs, and Mathematical Tables. New York, Dover, 1972.
\bibitem{bonami} A. Bonami, and J.-L. Clerc, Sommes de Ces\`aro et multiplicateurs des d\'eveloppement en harmoniques sph\'eriques, {\sl Trans. Amer. Math. Soc.} {\bf 183} (1973), 223--263.
\bibitem{bos} L. Bos, N. Levenberg, P. Milman, and B. Taylor, Tangential Markov inequalities on real algebraic varieties, {\sl Indiana Univ. Math. J.} {\bf 47} (1998) 1257-1272.
\bibitem{ditzian} Z. Ditzian, Fractional derivatives and best approximation, {\sl Acta Mathematica Hungarica} {\bf 81} (1998), 323--348.
\bibitem{koornwinder} T. Koornwinder, The addition formula for Jacobi Polynomials II. The Laplace Type Integral
Representation and the Product Formula, Math. Centrum Amsterdam,
Afd. Toegepaste Wisk. Rep. TV {\bf 133} (1972),\\ http://staff.science.uva.nl/~thk/art/.
\bibitem{muller} C. M\"uller, Spherical Harmonics, Springer-Verlag, 1966.
\bibitem{tao} A. Sikora and T. Tao,  Bochner-Riesz summability for analytic functions on the $m$-complex sphere, {\sl Communications in Analysis and Geometry} {\bf 12} (2004), 43--57.
%\bibitem{stein} E. M. Stein and G. Weiss, Introduction to Fourier Analysis on Euclidean Spaces, Princeton University Press, 1971.
\bibitem{szego} G. Szeg\"o, {\sl Orthogonal Polynomials}, Amer. Math. Soc., Providence, 1991.
\bibitem{wang} H.C. Wang, Two-point homogeneous spaces, Annals of Math. 55
(1952) 177-191.



%%%%%%%%%%%%%%%%%%%%%%%%%%%%%%%%%%%%%%%%%%%%%%%%%%%%%%%%%%%%%%%%%%%%%%
 \end{thebibliography}

%% not want to use model1c-num-names.bst.


\end{document}